\documentclass[12pt,leqno]{article}
\usepackage{graphicx}
\input amssym.tex
\begin{document}
\def\diamond{\square}
\def\pn{{\Bbb C}P(N)}
\def\pnh{\check{{\Bbb C}P(N)}}
\def\Sf{{\cal F}(2,d)}   	 
\def\C{\Bbb C}
\def\F{{\cal F}}	  		
\def\pr{{\bf Proof: }}
\def\S{\overline{\C}}
\def\P{{\cal P}}
\def\RR{\{F=0\}\cap\{G=0\}}
\def\INT{pGdF-qFdG}
\def\integ{\frac{F^p}{G^q}}
\def\integb{\frac{\bar{F}^p}{\bar{G}^q}}
\def\PP{{\cal P}_a \times {\cal P}_b}
\def\t{\tilde}
\def\In{{\cal I}(n,a,b)}
\def\Inm{{\cal I} _m(n,a,b)} 
\def\pro{{\Bbb C}P(n)}
\def\tpro{\tilde{{\Bbb C}P(n)}}
\def\Z{{\Bbb Z}}
\def\N{{\Bbb N}}
\def\H{{\cal H}}
\def\R{{\cal R}}
\def\integr{\frac{\t{F}^\frac{1}{q}}{\t{G}^\frac{1}{p}}}

\newtheorem{theo}{Theorem}[subsection]
\newtheorem{Theo}{Theorem}
\newtheorem{defi}{Definition}[subsection]
\newtheorem{coro}{Corollary}[subsection]
\newtheorem{lem}{Lemma}[subsection]
\newtheorem{claa}{Proposition}[section]
\newtheorem{exa}{Example}[subsection]
\newtheorem{rem}{Remark}[section]
\begin{center}
{\LARGE\bf On the Topology of Foliations with a First Integral\\}
\vspace{.25in}
{\large {\sc Hossein Movasati}\footnote{Partially supported by CNPq-Brazil
\\
\\
Math. classification: 14D99 - 57R30}
}
\end{center}
\begin{abstract}
The main objective of this article is to study the topology of the fibers of a generic rational function of the type $\integ$ in the projective space of dimension two. We will
prove that the action of the monodromy group on a single Lefschetz vanishing cycle $\delta$ generates the first homology group of a generic fiber of $\integ$.
In particular, we will prove that for any two Lefschetz vanishing cycles $\delta_0$ and $\delta_1$ in 
a regular compact fiber of $\integ$, there exists a monodromy $h$ such that
$h(\delta_0)=\pm \delta_1$. 
\end{abstract}
\setcounter{section}{-1}
\section{Introduction}
\ \ \ \ \
Let $F$ and $G$ be two homogeneous polynomials in $\C^{n+1}$. The following
function  is well-defined
\[
f=\integ :\pro\backslash\R \rightarrow \S
\]
\[
f(x)=
\frac{F(x)^p}{G(x)^q},\ x=[x_0;x_1;\cdots ;x_n]
\]
where $\R=\{F=0\}\cap\{G=0\}$, $\frac{deg(F)}{deg(G)}=\frac{q}{p}$ and $p$ and
$q$ are relatively prime numbers. We can view the fibration of $f$ as a
codimension one foliation in $\pro$ given by the 1-form
\[
\omega=\INT
\]
Let $\P_a$ denote the set of homogeneous polynomials of degree $a$ in $\C^n$. 
\begin{claa}
\label{gen}
There exists an open dense subset $U$ of $\P_a\times\P_b$ such that for any $(F,G)\in U$ we have:
\begin{enumerate}
\item
$\{F=0\}$ and $\{G=0\}$ are smooth varieties in $\pro$ and intersect each other 
transversally;
\item
The restriction of $f$ to ${\pro\backslash (\{F=0\}\cup \{G=0\})}$ has 
nondegenerate critical points, namely $p_1,p_2,\ldots,p_r$, 
with distinct images in $\S$, namely
$c_1,c_2,\ldots,c_r$ respectively.
\end{enumerate}
\end{claa}
Throughout the text the elements of $U$ will be called generic elements. We will prove this proposition in Appendix ~\ref{secgeneric}. Put 
\[
C=\{c_1,c_2,\ldots,c_r\}
\]
From now on, we will work with the function $f$ which has the generic
 properties as in
Proposition ~\ref{gen}. The foliation $\F$ associated to $f$ has the following
singular set
\[
Sing(\F)=\{p_1,p_2,\ldots,p_r,\R\}
\]
The value $0$ (resp. $\infty$) is a critical value of $f$ if and only if $p>1$ (resp. $q>1$). Let $A$ be a subset of $\{0,\infty\}$ which consists of only critical values.
For example if $p=1,q=1$ then $A$ is empty. The set of critical values of $\integ$
is $C\cup A$.  
\begin{claa}
$f$ is a $C^\infty$ fiber bundle map over $\S\backslash (C\cup A)$.
\end{claa}
We will prove this proposition in Section ~\ref{top}.
\\

The above proposition enables us to use the arguments of Picard-Lefschetz Theory
to study the topology of the fibers of $f$. But, for example, the critical fiber
$\{F=0\}$, when $p>1$, is not considered in that theory (as far as I know). To
overcome with this obstacle, we will construct a ramification map 
$\tau:\tpro\rightarrow\pro$ for the multivalued function $f^{\frac{1}{pq}}$. The
pull-back function $\t{f}^{\frac{1}{pq}}=f^{\frac{1}{pq}}\circ\tau$ of 
$f^{\frac{1}{pq}}$ is univalued and has no more the critical values $0$
and $\infty$. Next, we will embed the complex manifold $\tpro$ in some $\pn$ in
such a way that the pull-back foliation $\F$ is obtained by the 
intersection of the
hyperplanes of a generic Lefschetz pencil with $\tpro$.
\\
The study of the topology of an algebraic variety by intersecting it with
hyperplanes of a pencil has been started systematicly by Lefschetz in his famous
article \cite{lef}. We will use the arguments of this area of mathematics, specially the articles \cite{lam},\cite{che}, to
understand the topology of the leaves of $\F$. Note that the leaves of $\F$ do
not contain the points of the set $\{F=0\}\cap\{G=0\}$.   
\\
In the first section we will construct such ramification map $\tau$, and in the
second section we will review Picard-Lefschetz Theory. In the third section
we will apply our results to the foliation $\F$.
\\
 I want to use this opportunity to express my thanks to my advisor
 Alcides Lins Neto. Discussions with him during this work have provided valuable information 
 and new ideas. Also, thanks go to Cesar Camacho and Paulo Sad for their
 interest and support. I also thank Eduardo Esteves for his comments on algebraic geometry.  The author is also grateful to the exceptional 
 scientific atmosphere in IMPA that made this work possible.

\section{Ramification Maps}
\newtheorem{cla}{Proposition}[subsection]
\label{ramsec}
\ \ \ \ \
In the first part of this section we will introduce ramification maps with a
 normal crossing divisor. In the second part we will use a method which gives
  us some examples of ramification maps. This method will be enough for
  our purpose.
The following isomorphism will be used frequently during this section:

{\bf Leray (or Thom-Gysin) Isomorphism:} If a closed submanifold $N$ has 
pure real codimension $c$ in $M$, then there is an isomorphism 
\[
\tau :H_{k-c}(N)\tilde{\rightarrow}H_k(M,M\backslash N)
\]
 holding for any k, with the convention that $H_s(N)=0$ for
  $s<0$. Roughly speaking, given $x\in H_{k-c}(N)$, its image by this 
  isomorphism is obtained by thickening a cycle representing $x$, 
  each point of it growing into a closed $c$-disk transverse to $N$ in $M$ 
  (see \cite{che} p. 537). 
  \\
  Let $N$ be a connected codimension one submanifold of the complex manifold 
  $M$. Write the long exact sequence of the pair $(M,M\backslash N)$ as 
  follows:
\begin{equation}
\label{longler}
\cdots\rightarrow H_2(M,M\backslash N)
\stackrel{\sigma}{\rightarrow}H_1(M\backslash N)
\stackrel{i}{\rightarrow} H_1(M)\rightarrow \cdots
\end{equation}
where $\sigma$ is the boundary operator and $i$ is induced by inclusion. 
Since $N$ has real codimension two in $M$,  $H_2(M,M\backslash N)$ 
($\simeq H_0(N)\simeq\Z$) is generated by the disk $\Delta$ transverse to 
$N$ at a point $y\in N$. By the above long exact sequence it follows that 
if a closed cycle $x$ in $M\backslash N$ is homologous to zero in $M$ then
 it is homologous to a multiple of $\sigma (\Delta)=\delta$ in $M\backslash D$. 
The cycle $\delta$ is called a simple loop around the point $y\in N$ in
 $M\backslash N$.  
\subsection{Normal Crossing Divisors} 
\ \ \ \ \ 
 The following well-known fact will be used frequently:
\begin{cla}
\label{cover}
Let $\tau :\t{M}\rightarrow M$ be a finite covering map of degree $p$. 
Then the following statements are true:
\begin{enumerate}
\item
 $\tau _* :\pi_1(\t{M})\rightarrow \pi_1(M)$ is one to one, where $\pi_1(M)$
 denotes the fundamental group of $M$;
\item
 If $\pi_1(M)$ is abelian then $\pi_1(\t{M})$ is also abelian and 
 $\pi_1(M)^p\subset \tau_*(\pi_1(\t{M}))$, where $\pi_1(M)^p=\{\gamma ^p\mid\gamma\in\pi_1(M) \}$.
\end{enumerate}
\end{cla}
\pr 
The proof of the first statement can be found in \cite{ste}. 
If $\pi_1(M)$ is abelian then  
\[
\tau_*(aba^{-1}b^{-1})=\tau_*(a)\tau_*(b)\tau_*(a)^{-1}\tau_*(b)^{-1}=1
\]
where $a,b\in\pi_1(\t{M})$. The map $\tau_*$ is one to one and so 
$aba^{-1}b^{-1}=1$ which implies that $ab=ba$.
\\
For any closed path $a\in\pi_1(M,x)$, its inverse image by $\tau$ is a 
union of closed paths $a_1,a_2,\ldots,a_k$ in $\t{M}$. 
Choose a point $y$ in $\t{M}$ and points $x_i$ in the path $a_i$, 
for $i=1,\ldots, k$ such that $\tau(x_i)=\tau(y)=x$ and put
 $b_i=A_i^{-1}a_iA_i$ and $b=b_1b_2\cdots b_k$, where $A_i$ is a 
 path in $\t{M}$ which connects $y$ to $x_i$. The image of $A_i$ by $\tau$ is a 
 closed path in $M$ and $\pi_1(M)$ is abelian, therefore
\[
\tau_*(b)=\Pi\tau_*(A_i)^{-1}\tau_*(a_i)\tau_*(A_i)=\Pi \tau_*{a_i}=a^p
\]
The last equality is true because $\tau$ has degree $p$ and the paths $a_i$'s are 
the inverse image of $a$.$\diamond$
\\

In what follows, if the considered group is abelian, we use the additive notations 
of groups; for example instead of $a^p$ we write $pa$.  
    
\begin{defi}\rm 
Let $M$ be a complex manifold of dimension $n$. 
By a reduced normal crossing divisor we mean a union of finitely many 
connected closed submanifolds, namely  $D_1,D_2,\ldots,D_s$ of $M$ and of 
codimension one, which intersect each other transversally i.e., 
for any point $a\in M$ there is a local coordinate 
$(x,y)\in \C^k\times \C^{n-k},\ k\leq s$ around $a$ such that in this coordinate 
$a=(0,0)$ and for any $j=1,\ldots ,k$, the component $D_{i_j}$ is given by 
$x_j=0$,
 where $\{i_1,\ldots,i_k\}=\{ i\mid a\in D_i\}$. We say this coordinate 
 normalizing coordinate of $D$ at $a$.
  We will denote a reduced normal crossing divisor by 
  \[
  D=\sum_1^s D_i
  \]
When $D$ has only one component i.e., $s=1$, then $D$ is called simple. 
\end{defi}
The following fact is a direct consequence of the definition.
 \begin{cla}
\label{subncd}
Let $D$ be a reduced normal crossing divisor in a complex manifold  $M$. 
For any  subset $I$ of $\{1,2,\ldots,s\}$, the set 
\[
M_I=\cap_{i\in I} D_i
\]
is a complex manifold of codimension $\# I$ in $M$ and 
\[
D_I=D\cap M_I=\sum_{i\not\in I} D_i\cap M_I
\]
is a reduced normal crossing divisor in $M_I$.   
\end{cla}

In what follows for a given function $\tau :\t{M}\rightarrow M$,
and for any subset $x$ of $M$ (resp. a meromorphic function $x$ on $M$) 
we denote by $\t{x}$ the set $\tau^{-1}(x)$ 
(resp. the meromorphic function $x\circ\tau$ on $\t{M}$).
  
\begin{defi}\rm
\label{ncd}
Let $M$ and $\t{M}$ be two compact complex manifolds of the same dimension,
 $D=\sum_1^sD_i$ be a reduced normal crossing divisor in $M$ and 
 $p_1,p_2,\ldots,p_s$ be positive integer numbers  greater than one.
 The holomorphic map $\tau:\t{M}\rightarrow M$ is called a ramification map 
 with divisor $D$ and ramification index  $p_i$ at $D_i$ if
\begin{enumerate}
\item
$\tau^{-1}(D)=\t{D}$ is a reduced normal crossing divisor;
\item
 For any point $\t{a}\in \t{M}$ and a normalizing coordinate
 $(x,y)\in\C^k\times \C^{n-k}$ around $a=\tau(\t{a})$, there is a normalizing 
 coordinate $(\t{x},\t{y})\in\C^k\times \C^{n-k}$ around $\t{a}$ such that 
 in these coordinates $a=(0,0)$, $\t{a}=(0,0)$ and $\tau$ is given by: 
\[
(\t{x},\t{y})\rightarrow (x,y)=((\t{x}_1^{p_{i_1}},\t{x}_2^{p_{i_2}},
\ldots,\t{x}_k^{p_{i_k}}),\t{y})
\]
\end{enumerate}
\end{defi}
From the definition, we can see that
\begin{enumerate}
\item
The critical points and values of $\tau$ are $\t{D}$ and $D$, respectively;
\item
$\tau\mid_{\t{M}\backslash \t{D}}$ is a finite covering map of some degree $p$;
\item
For any subset $I$ of $\{1,2,\ldots,s\}$, if $M_I$ is not empty then 
$\Pi_{i\in I}p_i$ divides $p$.
\end{enumerate}
We have also the subramification maps of $\tau$, which is stated in bellow.
\begin{cla}
\label{subncd2}
Keeping the notations in Definition ~\ref{ncd} and Proposition 
~\ref{subncd}, for any ramification map $\tau:\t{M}\rightarrow M$ and $I$ 
a subset of $\{1,2,\ldots,s\}$ , the restriction of $\tau$ to $\t{M}_I$,
 namely $\tau_I$, is a ramification map with divisor $D_I$ and ramification 
 multiplicity $p_i$ at $D_i\cap M_I$, where $i\not\in I$. 
 Moreover, If $\tau$ is of degree $p$ then $\tau_I$ is of degree 
 $\frac{p}{\Pi_ {i\in I}p_i}$. 
\end{cla}

\begin{cla}
\label{cover2}
Let $\tau :\t{M}\rightarrow M$ be a ramification map of degree $p$ and with  
reduced normal crossing divisor $D$ and let $\pi_1(M\backslash D)$ be abelian.
 Then the following statement are true: 
\begin{enumerate}
\item
$\pi_1(\t{M})$ is abelian; 
\item
$p\pi_1(M)\subset \tau_*(\pi_1(\t{M}))$;
\item
 If $D$ is simple then $\tau_*: \pi_1(\t{M})\rightarrow \pi_1(M)$ 
is one to one;
\item
 If $D^*=\sum_{i=1}^{s^*}D^*_i$ is a reduced divisor in $M$ such that 
$D^*+D$ is a normal crossing divisor, then $\t{D^*}+\t{D}$ is also a 
normal 
crossing divisor. 
\end{enumerate}
\end{cla}

\pr The set $D$ is a finite union of some submanifolds of $M$ with real 
codimension greater than two, therefore every path in $\pi_1(M, x)$, \
where $x\in M\backslash D$, is homotopic to some path in 
$\pi_1(M\backslash D,x)$. Now the first and second statements are the direct 
consequences of  Proposition  ~\ref{cover}, Definition ~\ref{ncd} and the 
mentioned fact.
\\ 
Let $\Sigma$ be a small disk transverse to $D$ at $y\in D$ and 
$x\in\Sigma \cap M\backslash D$. 
Let also $a\in \pi_1(\t{M},\t{x})$, where $\tau(\t{x})=x$, and 
$\tau_*(a)$ be homotopic to zero in $M$ and $a$ do not intersect $\t{D}$. 
\\
Considering the long exact sequence ~(\ref{longler}) and the fact that 
$i(\tau_*(a))=0$, we conclude that $\tau_*(a)$ is homotopic to $k\delta$ 
in $M\backslash D$, where $k$ is an integer number and $\delta$ is a simple 
loop in $\Sigma$ around $y$. By Proposition ~\ref{cover}, this means 
that $a$ is homotopic to a closed path around $\t{y}$ in $\tau^{-1}(\Sigma)$, 
where $\tau(\t{y})=y$, which means that $a$ is homotopic to the point 
$\t{y}$ in $\t{M}$, and this proves the third statement.
\\
Let $a\in (\cap_{i\in I}D_i)\cap (\cap_{i\in I^*} D^*_{i})$, where
 $I\subset \{1,2,\ldots ,s\}$ and $I^*\subset \{1,2,\ldots ,s^*\}$.
  Choose a normalizing coordinate $(x,x^*,y)\in \C^r\times\C^{r^*}\times
  \C^{n-r-r^*}$
 around $a$ such that the components of $D$ (resp. $D^*$) through $p$ are 
 represented 
 by the coordinate $x$ (resp. $x^*$), where $r=\#I$ and $r^*=\#I^*$.
 Now the fourth statement is a direct consequence of  Definition ~\ref{ncd}.$\diamond$   

\subsection{ Construction of Ramification Maps}  
\ \ \ \ \
For any abelian group $G$ and a positive integer number $p$ define 
$G_p=G/pG$. We have the following properties with respect to $G_p$:
\begin{cla}
\label{group}
The following statements are true:
\begin{enumerate}
\item
Every morphism $f: G\rightarrow G'$ of abelian groups induces a natural 
morphism $f_p :G_p\rightarrow G'_p$;
\item
$(G_p)_q=G_{(p,q)}$, where $(p,q)$ denotes the greatest common divisor of 
$p$ and $q$;
\item
If $f:G\rightarrow G'$ is surjective then $f_p$ is also surjective. 
If $f$ is one to one then $f_p$ may not be one to one and so we cannot 
rewrite exact sequences of abelian groups by this change of groups and maps;
\item
Let $f:G\rightarrow G'$ be a morphism of abelian groups and $p$,$q$ be two
positive integer numbers. If $f$ is one to one, $pG'\subset f(G)$ and
 $(p,q)=1$ then $f_q$ is an isomorphism between $G_q$ and $G'_q$;
\end{enumerate}
\end{cla}

\pr We only prove the fourth statement, since others are trivial. There exist
 integer numbers $x,y$ such that $px+qy=1$. 
 \\
For any $a\in G'$ we have 
\[
a-q(ay)=p(ax)\in f(G)
\]
 and so $f_q$ is surjective. 
 \\
If $f_q(a)=0$ then $f(a)=qb$ for some $b\in G'$. 
We have 
\[
b=p(bx)+f(ay)=f(s),\ s\in G
\]
 which implies that $f(a-qs)=0$. The morphism $f$ is one to one and so 
 we have $a=qs$, which means that $f_q$ is one to one.$\diamond$
\\

The following statement gives us an example of ramification map with simple 
divisor.

\begin{cla}
Let $M$ be a complex manifold with $\pi_1(M)=0$ and $D$ be a simple divisor
 whose complement in $M$ has abelian fundamental group. 
 Let also $p'$ be a positive integer number and 
 \[
 \pi_1(M\backslash D)_{p'}=\Z_p
 \]
Then there exists a degree $p$ ramification map with divisor $D$ 
and ramification multiplicity $p$ at $D$.
\label{newman}
\end{cla}

\pr  For any $e\in D$ define $\pi_1(M\backslash D,e)=\{0_e\}$ and
\[
\tilde{M}=\cup_{e\in M} \pi_1(M\backslash D,e)_{p'}
\]
$\tilde{M}$ has the structure of a complex manifold. 
For any $[a_e]\in \pi_1(M\backslash D,e)$ we must define a base open set 
and a chart map around $[a_e]$. Consider two cases:
\\
1. $e\in M\backslash D$
\\
Let $V_e$ be a simply connected open neighborhood of $e$ in 
$M\backslash D$. The following function is well-defined:
\[
\eta :V_e\rightarrow \t{M}
\]
\[
\eta(y)=[A_{ey}a_eA_{ey}^{-1}]
\]
where $A_{ey}$ is a path which connects $e$ to $y$ in $V_e$. The image of $\eta$ is
 a base open set around $a_e$ and $\eta$ is a chart map.
\\
2. $e\in D$
\\
 Let $(V_e,(x,y))$, $(x,y)\in (\C^{n-1}\times \C ,0)$, be a coordinate around
  $e$ such that in this coordinate $e=(0,0)$ and $D$ is given by $y=0$. 
  By Leray isomorphism, for any $e'=(x',y')\in V_e$ the group 
  $\pi_1(M\backslash D ,e')_{p'}\simeq \Z_p$ is generated by a simple 
  loop around $(x',0)$ in 
  \[
  \Sigma _{(x',0)}=\{(x,y)\in (\C^{n-1}\times \C, 0)\mid x=x'\}
  \]
 In particular, we have  $p\mid p'$.
   This gives us the following construction of a chart map 
  around $[a_e]$:
\\
 For any $y_0\in (\C,0)$, let $j(y_0)$ be a point in $(\C,0)$ such that
 \[
 y_0^p=j(y_0)^p\ \&\ 0\leq arg(j(y_0))< \frac{2\pi}{p}
 \]
 and let $\delta_{y_0}$ be the path which connects $y_0$ to $j(y_0)$ in 
 $\{y\in (\C,0)\mid y^p=y_0^p\}$ in the clock direction.
 \\
 The image of $\delta_{y_0}$ by the map $i(y)=y^p$, $(\delta_{y_0})^p$, 
 is a closed path with initial and end point $y_0^p$ and so the following
 function is well-defined
\[
\eta :(\C^{n-1}\times\C ,0)\rightarrow \t{M}
\]
\[
\eta (x,y)=\{x\}\times (\delta_y)^p
\]
The image of $\eta$ is a base open set around the point $e$ and $\eta$ is a 
chart map. The reader can verify easily that $\t{M}$ with these base open sets
and chart maps is a complex manifold and
the natural function $\tau :\t{M}\rightarrow M$ 
is the desired ramification map.
$\diamond$       

\begin{theo}
\label{ramtheo}
Let $M$ be a complex manifold with $\pi_1(M)=0$ and $D=\sum_{i=1}^sD_i$ be a reduced normal 
crossing divisor such that  the complement of each $D_i$ in $M$ has
 abelian fundamental group. 
 Let also $p'_1,p_2',\ldots p_s'$ be  positive integer numbers which 
 are prime to each other. Put $p_i=\# \pi_1(M\backslash D_i)_{p_i'}$.
 Then there is a degree $p_1p_2\cdots p_s$ ramification map with divisor $D$
  and ramification multiplicity $p_i$ at $D_i$, $i=1,2,\ldots,s$.
\end{theo}

\pr The proof is by induction on $s$. For $s=1$ it is Proposition 
~\ref{newman}. 
\\
Suppose that the theorem is true for $s-1$. 
Let $\tau :\t{M}\rightarrow M$ be a degree $p_1$ ramification map  
with simple divisor $D_1$ and multiplicity $p_1$ at $D_1$. 
We check the assumptions of the theorem for the divisor 
$\t{D_2}+\ldots+\t{D_s}$ in  the manifold $\t{M}$,
 to apply the hypothesis of the induction.
 \\
By the third part of Proposition ~\ref{cover2}
$\tau_* : \pi_1(\t{M})\rightarrow \pi_1(M)$ is one to one and by hypothesis
 $\pi_1(M)=0$, therefore $\pi_1(\t{M})=0$.
\\
Applying  Proposition ~\ref{cover2} to the ramification map 
$\tau\mid_{\t{M}\backslash \t{D_i}}$, we see 
that $\pi_1(\t{M}\backslash \t{D_i})$ is abelian;
 also  4 of Proposition ~\ref{cover2}  
 implies that $\t{D}$ is a normal crossing divisor. 
\\
The morphism 
\[
\tau _*: \pi_1(\t{M}\backslash \t{D_i})\rightarrow \pi_1(M\backslash D_i)
\]
 is one to one and by 2 of Proposition ~\ref{cover2} 
 \[
 p_1\pi_1(M\backslash D_i)\subset \tau_* (\pi_1(\t{M}\backslash \t{D_i})
 \]
 But $g.c.d.(p_1,p'_i)=1$ and so by 4 of Proposition ~\ref{group}
 we have
\[
     \pi_1(\t{M}\backslash \t{D_i})_{p_i'}\simeq 
     \pi_1(M\backslash D_i)_{p_i'}
\] 
Now we can apply the hypothesis of the induction to $\t{M}$ and 
$D'=\sum_2^sD_i$. There exists a degree $p_2\cdots p_s$ ramification map 
$\tau':\t{M}'\rightarrow \t{M}$ with divisor $D'$ and multiplicity $p_i$ at 
$D_i$, $i=2,\ldots ,s$. 
The reader can check that the  map $\tau\circ\tau '$ is the desired 
ramification map. $\diamond$
\subsection{Multivalued Functions}
\ \ \ \ \
To study multivalued functions, we will need to study a certain class of
ramification maps. 
First, we give the precise definition of multivalued functions.
\begin{defi}\rm
Let $\tau:\t{M}\rightarrow M$ be a degree $p$ holomorphic map between two 
complex manifold $\t{M}$ and $M$ i.e., $\tau$ is a finite covering map of 
degree $p$ out of its critical points.
 Every meromorphic function $g$ on $\t{M}$ is called a $p$-valued meromorphic
 function on $M$.
  Roughly speaking, the image of a point $x\in M$ under $g$ is the
   set $g(\tau^{-1}(x))$. The map $\tau$ is called the ramification map of $g$ 
   and the set of critical values of the map $\tau$ is called the ramification 
   divisor of the multivalued function $g$.   
\end{defi}

Given a complex manifold $M$, a meromorphic function $f$ on it and 
an integer number $N$. 
Can we  construct a ramification map of the multivalued function 
$f^{\frac{1}{N}}$ according to the above definition?
 Here we will answer to this question for some limited classes of meromorphic 
 functions.

Let $f$ be a meromorphic function on the manifold $M$ 
and $\tau:\t{M}\rightarrow M$ be a ramification map with reduced normal 
crossing divisor $D=\sum_{i=1}^sD_i$ and multiplicity $p_i$ at $D_i$, 
$i=1,2,\ldots,s$.
 Let also $div(f)=\sum m_j V_j$. Then
\[
div(\t{f})=\sum a_jm_jV_j 
\]
where $\tilde{f}=f\circ\tau$, $a_j=p_i$ if $V_j=D_i$ for some $i$ and $a_j=1$ otherwise.
\begin{cla}
Keeping the notations used above, suppose that $H^1(\t{M} ,\Z _N)=0$ and 
$N\mid div(\t{f})$ i.e., $N$ divides the multiplicities of the components of 
$div(\t{f})$, where $N$ is a positive integer number. 
Then $\t{f}^{\frac{1}{N}}$ is a well-defined meromorphic function on $\t{M}$.
Therefore, we can view $f^{\frac{1}{N}}$ as a multivalued function on $M$ with
the ramification map $\tau$.
\end{cla}

\pr For any point $x\in \t{M}$ there is a neighborhood $V_i$ of $x$
 such that in this neighborhood $\t{f}=g_i^N$, where $g_i$ is a meromorphic 
 function on $V_i$. Let $c_{ij}=\frac{g_i}{g_j}$, then $c_{ij}^N=1$. 
 Since $H^1(\t{M} ,\Z _N)=0$, there exist complex numbers $c_i$'s  
 such that $c_{ij}=\frac{c_i}{c_j}$. 
 Now $g\mid_{V_i}=\frac{g}{c_i}$'s define the global meromorphic function which is the 
 desired candidate for $\t{f}^{\frac{1}{N}}$. 
 $\diamond$    
\section{Picard-Lefschetz Theory}
\label{picsec}
\ \ \ \ \
In 1924 S. Lefschetz published his famous article \cite{lef} on the topology of 
algebraic varieties. In his article, in order to study the topology of an
algebraic variety, he considered a pencil of hyperplanes in general position 
with respect to that variety. 
Many of the Lefschetz intuitive arguments are made precise by appearance 
a critical fiber bundle map. In the first part of this section we introduce
the basic concepts of Picard-Lefschetz Theory and in the second part we 
introduce the Lefschetz pencil and state our two basic theorems
~\ref{ggenerator},~\ref{monodromy}. This section is mainly based on the articles \cite{lam},\cite{che}. Homologies  are considered in an
arbitrary field of characteristic zero except it mentioned explicitly.
\subsection{Critical Fiber Bundle Maps} 
\ \ \ \ \
The following theorem gives us a huge number of fiber bundle  maps. 
\begin{Theo}
(Ehresmann's Fibration Theorem \cite{ehr}). 
Let $f:Y\rightarrow B$ be a proper submersion between the manifolds $Y$ and $B$.
 Then $f$ fibers $Y$ locally trivially i.e., for every point $b\in B$ 
 there is a neighborhood $U$ of $b$ and a $C^\infty$-diffeomorphism 
 $\phi :U\times f^{-1}(b) \rightarrow f^{-1}(U)$ such that 
 $f\circ\phi =\pi _1 =$ the first projection. Moreover if $N\subset Y$ 
 is a closed submanifold such that $f\mid _N$ is still a submersion then 
 $f$ fibers $Y$ locally trivially over N i.e., the diffeomorphism  
 $\phi$ above  can be chosen to carry  $U\times (f^{-1}(b)\cap N)$ 
 onto $f^{-1}(U)\cap N$.
\end{Theo}

The map $\phi$ is called the fiber bundle trivialization map. 
Ehresmann's theorem can be rewrite for manifolds with boundary and also 
for stratified analytic sets. In the last case the result is known as the 
Thom-Mather theorem.

In the above theorem let $f$ not be submersion, and let $C'$ be the union of 
critical values of $f$ and critical values of $f\mid _N$, and $C$ be the 
closure of $C'$ in $B$. By a critical point of the map $f$ we mean the point 
in which $f$ is not submersion. Now we can apply the theorem to the function 
\[
f: Y\backslash f^{-1}(C)\rightarrow B\backslash C=B'
\]
For any set $K\subset B$, we use the following notations
\[
Y_{K}=f^{-1}(K) ,\ Y'_K=Y_K\cap N,\ L_K=Y_K\backslash Y'_K
\]
and for any point $c\in B$, by $Y_c$ we mean the set $Y_{\{c\}}$. By  
\[
f:(Y,N)\rightarrow B
\]
 we mean the mentioned map and we call it the critical fiber bundle map. 

\begin{defi}\rm
Let $A\subset R\subset S$ be topological spaces. $R$ is called a strong 
deformation retract of $S$ over $A$ if there is a continuous map 
$r:[0,1]\times S\rightarrow S$ such that 
\begin{enumerate}
\item
$r(0,.)=id$;
\item
$r(1,x)\in R\ \&\ r(1,y)=y \ \forall x,y\in S,\ y\in R$; 
\item
$r(t,x)=x \ \forall t\in[0,1],\ x\in A$.
\end{enumerate}
Here $r$ is called the contraction map. 
In a similar way we can do this definition for the pairs of spaces 
$(R_1, R_2)\subset (S_1,S_2)$, where $R_2\subset R_1$ and $S_2\subset S_1$.
\end{defi}

We use the following important theorem to define generalized vanishing cycle 
and also to find  relations between the homology groups of $Y\backslash N$ 
and the generic fiber $L_c$ of $f$.
\begin{theo}
\label{contraction}
Let $f:Y\rightarrow B$ and $C'$ as before, $A\subset R\subset S\subset B$ and $S\cap C$ be a subset of the interior 
of $A$ in $S$, then every retraction from $S$ to $R$ over $A$ can be lifted
 to a retraction from $L_S$ to $L_R$ over $L_A$. 
\end{theo}

\pr According to Ehresmann's fibration theorem 
$f:L_{S\backslash C}\rightarrow S\backslash C$ is a $C^\infty$ locally trivial 
fiber bundle. The homotopy covering theorem, see 14,11.3\cite{ste}, implies 
that the contraction of $S\backslash C$ to $R\backslash C$ over 
$A\backslash C$ can be lifted so that $L_{R\backslash C}$ becomes a strong 
deformation retract of  $L_{S\backslash C}$ over $L_{A\backslash C}$. 
Since $C\cap S$ is a subset of the interior of $A$ in $S$,
 the singular fibers can be filled in such a way that 
 $L_R$ is a deformation retract of $L_S$ over $L_A$.$\diamond$
 \\
 
 {\bf Monodromy:} Let $\lambda$ be a path in $B'=B\backslash C$ with the initial 
and end points $b_0$ and $b_1$. In the sequel by $\lambda$ we will mean both 
the path $\lambda :[0,1]\rightarrow B$ and  the image of $\lambda$; 
the meaning being clear from the text.
\begin{cla}
\label{isotopy}
There is an isotopy 
\[
H:L_{b_0}\times[0,1]\rightarrow L_{\lambda}
\]
such that for all $x\in L_{b_0},\ t\in [0,1]$ and  $y\in N$
\begin{equation}
\label{iso}
 H(x,0)=x,\ H(x,t)\in L_{\lambda (t)},\  H(y,t)\in N
\end{equation}
For every $t\in [0,1]$ the map $h_t=H(.,t)$ is a homeomorphism between 
$L_{b_0}$ and $L_{\lambda (t)}$. The different choices of $H$
 and paths homotopic
 to $\lambda$ would give the class of homotopic maps 
 \[
 \{h_{\lambda} : L_{b_0} \rightarrow L_{b_1}\} 
 \]
 where $h_\lambda=H(.,1)$.
 \end{cla}

\pr The interval $[0,1]$ is compact and the local trivializations of
 $L_\lambda$ can be fitted together along $\gamma$ to yield an isotopy $H$.
  $\diamond$

The class $\{h_{\lambda} : L_{b_0} \rightarrow L_{b_1}\} $ defines the maps
\[
h_{\lambda}:\pi _*(L_{b_0})\rightarrow\pi_*(L_{b_1})
\]
\[ 
h _{\lambda} :H_*(L_{b_0}) \rightarrow H_*(L_{b_1})
\]
In what follows we will consider the homology class of cycles, but many of 
the arguments can be rewritten for their homotopy class. 
\begin{defi}\rm
For any regular value $b$ of $f$, we can define 
\[
h:\pi _1(B',b)\times H_*(L_b)\rightarrow H_*(L_{b})
\]
\[
h(\lambda ,.)=h_\lambda(.)
\]
$\pi_1(B',b)$ is called the monodromy group and its action 
$h$ on $H_*(L_b)$ is called the action of monodromy on the homology groups of
$L_b$.
\end{defi}
Following the article \cite{che}, we give the generalized definition of vanishing cycles.      
\begin{defi}
Let $K$ be a subset of $B$ and  $b$ be a point in $ K\backslash C$. 
Any relative k-cycle of $L_K$ modulo $L_b$ is called  a k-thimble above 
$(K,b)$ and its boundary in $L_b$ is called a vanishing $(k-1)$-cycle above $K$.
\end{defi}

\begin{figure}[t]
\begin{center}
\includegraphics{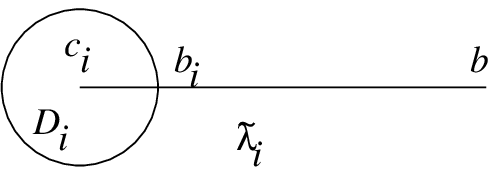}
\caption{}
\label{alc}
\end{center}
\end{figure}

Let us consider the case that we will need. Let $Y$ be a complex  compact manifold,
 $N$ be a submanifold of $Y$ of codimension one,
  $B=\S$ and $f$ be a holomorphic function.
 The set of critical values of $f$, $C$, is a finite set.
 \\ 
Let $c_i\in C$ (which is an 
isolated point of $C$ in $\S$), $D_i$ be an small 
disk around $c_i$ and $\tilde{\lambda _i}$ be a path in $B'$ which connects
 $b\in B'$ to $b_i\in \partial D_i$. 
 Put $\lambda _i$ the path $\tilde{\lambda _i}$ plus the path 
 which connects $b_i$ to $c_i$ in $D_i$ (see Figure ~\ref{alc}). Define the set $K$ in the three ways 
 as follows:

\begin{equation}
K^s= 
	 \left\{ \begin{array}{lr}
		\lambda_i & s=1\\
		\lambda_i \cup D_i & s=2 \\
		\tilde{\lambda}_i\cup \partial D_i & s=3\\
		       \end{array} 
	\right.
\end{equation}
In each case we can define the vanishing cycle in $L_b$ above $K^s$.
$K^1$ and $K^3$ are subsets of $K^2$ and so the vanishing cycle above $K^1$ and
$K^3$ is also vanishing above $K^2$. 
In the first case we have the intuitional concept of vanishing cycle. 
If $c_i$ is a critical point of $f\mid_N$ we can see that the vanishing cycle
above $K^2$ may not be vanishing above $K^1$.  
\\
The third case gives us the vanishing cycles obtained by a monodromy around 
$c_i$. In this case we have the Wang isomorphism
\[
v: H_{k-1}(L_b)\tilde{\rightarrow} H_k(L_K,L_b)
\]
Roughly speaking, The image of the cycle $\alpha$ by $v$ is the footprint
of $\alpha$, taking the monodromy around $c_i$. 
Let $\gamma_i$ be the closed path which parameterize $K_3$ i.e., 
$\gamma _i$ starts from $b$, goes along $\tilde{\lambda}_i$ until $b_i$, 
turn around $c_i$ on $\partial D_i$ and finally comes back to $b$ along 
$\tilde{\lambda_i}$. 
Let also $h_{\gamma _i}:H_k(L_{b_i}) \rightarrow H_k(L_{b_i})$ 
be the monodromy around the critical value $c_i$. 
It is easy to check that 
\[
\sigma\circ v= h_{\lambda_i}-I
\]
 where $\sigma$ is the boundary operator, 
therefore the cycle $\alpha$ is a vanishing cycle above $K^3$ if and 
only if it is in the image of $h_{\lambda_i}-I$.
For more information about the generalized vanishing cycle the reader is
referred to \cite{che}.
\\

{\bf Lefschetz Vanishing Cycle:} Let $f$ have a nondegenerate critical point
 $p_i$ in $Y\backslash N$ and $p_i$ be the unique critical point of 
 $f:(Y,N)\rightarrow \S$ within $Y_{c_i}$, where $c_i=f(p_i)$.
 
\begin{cla}
\label{uvan}
In the above situation, the following statements are true:
\begin{enumerate}
\item
 For all $k\not = n$ we have $H_k(L_{\lambda_i} ,L_b)=0$. This means that there 
 is no $(k-1)$-vanishing cycle along $\lambda_i$ for $k\not =n$;
\item
$H_n(L_{\lambda_i} ,L_b)$ is infinite cyclic generated by a hemispherical 
homology class $[\Delta_i]$ which is called the Lefschetz thimble and its 
boundary is called the Lefschetz vanishing cycle;   
\item
 Let $\lambda_i'$ be another path which connects $b$ to $c_i$ in $B'$ and is
  homotopic to $\lambda_i$ in $B'$ (with fixed initial and end point), 
then we have the same,  up to homotopy and
  change of sign, Lefschetz vanishing cycle in $L_b$.   
\end{enumerate}
\end{cla}
 For the proof of above Proposition see 5.4.1 of \cite{lam}  .
\\
By a hemispherical homology class, we mean the image of a generator of 
infinite cyclic group $H_n({\Bbb B}^n, {\Bbb S}^{n-1})$ under the homeomorphism 
induced by a continuous mapping of the closed $n$-ball ${\Bbb B}^n$ into $L_{\lambda_i}$ which sends its boundary, the $(n-1)$-sphere ${\Bbb S}^{n-1}$, to $L_{b}$. Let $B$ be a
small ball around $p_i$ such that in $B$ we can write $f$ in the Morse form
\[
f=c_i+x_1^2+x_2^2+\cdots+x_n^2
\]
For $b$ such that $b-c_i$ is positive real, the Lefschetz vanishing cycle 
in the fiber $L_b$ is given by:
\[
\delta_i =\{ (x_1,\cdots ,x_n)\in {\Bbb R}^n\mid \sum x_j^2=b-c_i\}
\]
which is the boundary of the thimble
\[
\Delta_i =\{(x_1,\cdots ,x_n)\in {\Bbb R}^n\mid \sum x_j^2\leq b-c_i\}
\]
In the above situation the monodromy $h_i$ around the critical value $c_i$ is given by the Picard-Lefschetz formula
\[
h(\delta)=\delta+(-1)^{\frac{n(n+1)}{2}}<\delta,\delta_i>\delta _i,\ \ \ \delta\in H_{n-1}(L_b)
\]
where  $<.,.>$ denotes the intersection number of two cycles in $L_{b}$.
\\
{\bf Remark:} In the above example vanishing above
 $K^1$ and $K^2$ are the same. Also by the Picard-Lefschetz formula the reader can verify that three types of the definition of a 
 vanishing cycle coincide.  
 In what follows by vanishing along the path $\lambda_i$ we will mean vanishing
 above $K^2$.
\subsection{Vanishing Cycles as Generators} 
\begin{figure}[t]
\begin{center}
\includegraphics{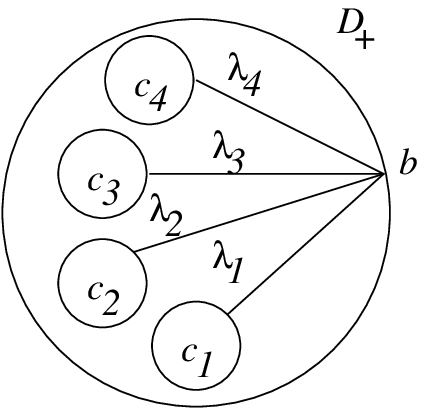}
\caption{}
\label{monono1}
\end{center}
\end{figure}

\ \ \ \ \
Now let $\{c_1,c_2,\ldots,c_s\}$ be a subset of the set $C$ of  critical values of
$f$, and $b\in \S\backslash C$. 
Consider a system of $s$  paths $\lambda _1,\ldots,\lambda _s$ starting 
from $b$ and ending at $c_1, c_2,\ldots , c_s$, respectively, and such that:
\begin{enumerate}
\item
 each path $\lambda _i$ has no self intersection points  ; 
\item
two distinct path $\lambda _i$ and $\lambda _j$ meet only at their 
 common origin $\lambda _i(0)=\lambda _j(0)=b$ (see Figure ~\ref{monono1}). 
\end{enumerate}
 This system of paths is called  a distinguished system of paths. 
 The set of vanishing cycles along the paths 
 $\lambda _i,\ i=1,\ldots ,s$ is called a distinguished set of vanishing 
 cycles related to the critical points $c_1,c_2,\ldots,c_s$.

\begin{theo}
\label{ggenerator}
Suppose that $H_{k-1}(L_{\S\backslash\{a\}})=0$ for some positive integer number $k$
 and $a\in \S$, which may be a critical value. 
 Then a distinguished set of vanishing $(k-1)$-cycles 
 related to the critical points in the set 
  $C\backslash \{a\}=\{c_1,c_2,\ldots,c_r\}$ generates $H_{k-1}(L_b)$.
\end{theo}

\pr We use the arguments of the article \cite{lam} Section 5. Note that in our case the fiber is $L_b=Y_b\backslash N$ and not $Y_b$.
\\
We consider our system of distinguished paths inside a large disk 
$D_+$ so that $a\in \S\backslash \overline{D}_+$, the point $b$ 
is in the  boundary of $D_+$ and all critical values $c_i$'s in 
$C\backslash \{a\}$ are interior points of $D_+$. 
Small disks $D_i$ with centers $c_i$ $i=1,\cdots ,r$ are chosen so that they
 are mutually disjoint and contained in $D_+$. 
 Put 
 \[
 K_i=\lambda _i\cup D_i,\ K=\cup_{i=1}^r K_i
 \]
   The pair $(K,b)$ is a strong deformation retract of $(D_+,b)$ 
   and so by Theorem ~\ref{contraction} $(L_K,L_b)$ is a strong deformation retract of
    $(L_{D_+},L_b)$. The set $\tilde{\lambda}=\cup \tilde{\lambda_i}$ 
    can be retract within itself to the point $b$ and so $(L_K,L_b)$ and 
    $(L_K,L_{\tilde{\lambda}})$ have the same homotopy type. 
    By the excision theorem (see \cite{mas}) we conclude that 
\[
H_k(L_{D_+},L_{b})\simeq\sum _{i=1}^{r} H_k(L_{K_i},L_{b})\simeq\sum _{i=1}^{r} H_k(L_{D_i},L_{b_i})
\]
Write the long exact sequence of the pair $(L_{D_+},L_b)$:
\begin{equation}
\ldots\rightarrow H_k(L_{D_+})
\rightarrow H_k(L_{D_+},L_b) \stackrel{\sigma}{\rightarrow } H_{k-1}(L_b) 
\rightarrow H_{k-1}(L_{D_+})\rightarrow \ldots
\label{long}
\end{equation}
Knowing this long exact sequence, 
it is enough to prove that $H_{k-1}(L_{D_+})=0$. 
A contraction from $\S\backslash \{a\}$ to $D_+$ can be lifted
 to the contraction of $L_{\S\backslash \{a\}}$ to $L_{D_+}$ 
 which means that $L_{D_+}$ and $L_{\S\backslash\{a\}}$ have the same 
 homotopy type and so by the  hypothesis $H_{k-1}(L_{D_+})=0$.$\diamond$ 
\subsection{Lefschetz Pencil}
\ \ \ \ \
In this section we repeat some notations and propositions  of \cite{lam} Section
2. All the proofs can be found there.  
\\   
The hyperplanes of $\pn$ are  points of the dual projective space $\pnh$. 
We use the following notation:
\[
H_y \subset \pn\ , \ y\in\pnh
\]
Let $X$ be a closed irreducible subvariety of $\pn$ and let $X_e \subset X$ 
be the nonempty open subset of its regular points. Define 
\[
     V_X^{'}=\{(x,y)\in \pn \times \pnh \mid x\in X_e\ \&\ H_y 
     \hbox{ is tangent to
     $X$ at $x$}\}
\]
This is a quasiprojective subset of $\pn \times \pnh$, because the set
\[
	\tilde{V}=\{(x,y)\in \pn \times \pnh \mid 
	\hbox{ $x$ is a singular point of $X$
        or $H_y$ is tangent to}
\]
\[
 	X \hbox{ at } x \}
\]
is closed in $\pn \times \pnh $ and $V_X^{'}$ is a Zariski open in 
$\tilde{V}$. The closure $V_X$ of $V_X^{'}$ is called the tangent hyperplane 
bundle of $X$. Consider the second projection 
\[
	\pi _2 : V_X \rightarrow \pnh \ ,\ (x,y)\rightarrow y
\]
its image $\check{X}$ is a closed irreducible subvariety of $\pnh$ of dimension 
at most $n-1$ which is called the dual variety of $X$.
If $X$ is a smooth variety then
\[
 \check{X}=\{y\in \pnh \mid H_y \hbox{ is tangent to $X$ at some point }\}
\]
In general $\check{X}$
 has singularities even if $X$ does not. 
 If $dim(\check{X})=N-1$ the degree of $\check{X}$ is well-defined and if 
 $dim(\check{X})<N-1$ we define $deg(\check{X})=0$.  
\begin{cla}
(Duality Theorem \cite{lam} 2.2) The tangent hyperplane bundles of $X$ and $\check{X}$ coincide
\[
V_X=V_{\check{X}} \hbox{ and hence } \check{\check{X}}=X
\]
\end{cla}

A pencil in $\pn$ consists of  all hyperplanes which contain a fixed 
$(N-2)$-dimensional projective space  $A$, which is called the axis 
of the pencil. We denote a pencil by 
\[
\{H_t\}_{t\in G}
\]
or $G$ itself, where $G$ is a projective line in $\pnh$. 
\\
The pencil $\{H_t\}_{t\in G}$ is in general position with respect to $X$
 if $G$ is in general position with respect to $\check{X}$. 
  From now on, fix a  pencil $\{H_t\}_{t\in G}$ 
  in general position with
  respect to $X$.
  
\begin{cla}
( \cite{lam}, 1.6.1) The axis A intersects X transversally.
\end{cla}
For the pencil $\{H_t\}_{t\in G}$ put 
\[
X_t=X\cap H_t,\ L_t=X_t\backslash A ,\ C=G\cap \check{X}=
\{c_1,c_2,\ldots,c_r\},\ X'=A\cap X  
\]
We will sometimes parameterize $G$ by $\S$ and denote the pencil by 
  $\{H_t\}_{t\in \S}$.
In order to have a map whose level surfaces are the $X_t$'s,  we need to do a "blow up"
 along the  variety $X'$. Let 
\[
Y=\{(x,t)\in X\times \S \mid x\in H_t\}
\]
There are two projections
\[
X\stackrel{p}{\leftarrow} Y \stackrel{f}{\rightarrow} \S
\]
Put $Y'=p^{-1}(X')=X'\times \S$ then
\begin{cla}
\label{fibration}
(\cite{lam} 1.6.2, 1.6.3, 1.6.4) If $X$ is a smooth variety then
\begin{enumerate}
\item
The modification $Y$ of $X$ is smooth and irreducible;
\item 
$p$ is an isomorphism between $Y\backslash Y'$ and $X\backslash X'$ and also 
an isomorphism between $f^{-1}(t)$ and $X_t$;
\item
For every critical value $c_i,\ i=1,\ldots ,r$  of $f$, the hyperplane 
 $H_{c_i}$ has a unique tangency of order two with $X$ which lies out of $A$.
 The other hyperplanes $H_c,\ c\not\in C$ are transverse to $X$;
\item
The projection $f:Y\rightarrow \S$ has $r=deg(\check{X})$ nondegenerate 
critical points $p_1,\ldots, p_r$ in $Y\backslash Y'$ such that $f(p_i)=c_i$'s 
are distinct 
values in $\S$.
\end{enumerate}
 \end{cla}

Now we have the critical fiber bundle map 
$f:(Y, Y')\rightarrow \S$. Note that $f\mid _{Y'}$ has no critical points. 
 We conclude that the natural function $f:X\backslash A\rightarrow \S$ is a
 fiber bundle map over $\S \backslash C$.

\begin{defi}\rm 
We can view $f:X\rightarrow \S$ as a meromorphic function on $X$. $f$ is called
the Lefschetz meromorphic function. The foliation induced by the pencil
$\{H_t\}_{t\in \S}$ is called the Lefschetz foliation.
\end{defi}

\begin{cla}
\label{equiv}
The pencil $\{H_t\}_{t\in G}$ is in general position with respect to $X$ if
 and only if
 \begin{enumerate}
 \item
 Choosing a good parameterization of $G$, $(f)_\infty=f^{-1}(\infty)$ and 
$(f)_0=f^{-1}(0)$ are smooth irreducible varieties and intersect each
 other transversally and ;
\item
$f$ has only nondegenerate critical points with distinct images.
\end{enumerate}
 \end{cla}

\pr If $\{H_t\}_{t\in G}$ is in general position with respect to $X$, then the axis
of the pencil, $A$, intersects $X$ transversally. Knowing that
\[
(f)_0=H_0\cap X,\ (f)_\infty=H_\infty\cap X, \ H_0\cap H_\infty=A
\]
we conclude that $(f)_0$ and $(f)_\infty$ intersect each other transversally.
The second statement is Proposition ~\ref{fibration}, 4.
\\
Now suppose that $f$ satisfies 1 and 2. Suppose that $G$ is not
in general position with respect to $\check{X}$ at $x\in G\cap\check{X}$.   
We can distinguish two cases
\\
1. $x$ is a smooth point of $\check{X}$ and $G$ is tangent to $\check{X}$ at
$x$;
\\
Let $H_s$ be a hyperplane which passes through $x$, contains $G$ and is tangent
to $\check{X}$ at $x$. By Duality Theorem the information
\[
x\in G\subset H_s, \ \hbox{ $H_s$ is tangent to $\check{X}$ at $x$}
\]
can be translated to 
\[
s\in A\subset H_x, \ \hbox{ $H_x$ is tangent to $X$ at $s$}
\]
But this contradicts the first statement.
\\
2. $x$ is a singular point of $\check{X}$;
\\
This case also cannot happen. By the argument used in the proof of 1.6.4 of \cite{lam}, we have: $x$ is a smooth point of $\check{X}$, if and only if, $H_x$ has a unique tangency point of order two with $X$.$\diamond$
\begin{theo}
Suppose that the pencil $\{H_t\}_{t\in G }$ is in general position with respect to $X$ 
and Let $a$ be a point in $\S\backslash C$. Then for every $b\in\S\backslash C$ \label{generator}
\begin{enumerate}
\item
 $H_i(L_b)\simeq H_i(X\backslash H_a),\ i\not = n, n-1$ 
  \item
If $H_{n-1}(X\backslash H_a)=0$, then a distinguished set of vanishing 
cycles related to the critical values $c_1,c_2,\ldots,c_r$ generates the 
group $H_{n-1}(L_b)$. 
\end{enumerate}
\end{theo}

\pr This is a direct consequence of Theorem ~\ref{ggenerator} and  
Theorem ~\ref{uvan} and the long exact sequence ~(\ref{long}).$\diamond$
\\

{\bf Blow up:} Fix the point $b\in\pro$. All lines through $b$ in $\pro$ 
form a projective space of dimension $n-1$, namely $P$. Define
\[
\t{\pro}=\{(x,y)\in \pn\times P\mid x\in y\}
\]
$\t{\pro}$ is a smooth subvariety of $\pn\times P$. We have two natural
 projections
\[
\pro\stackrel{i}{\leftarrow} \t{\pro} \stackrel{f}{\rightarrow} P
\]
The reader can check that $f$ is an isomorphism between $i^{-1}(b)$ and $P$ 
and $i$ is an isomorphism between 
$\pro\backslash\{b\}$ and $\t{\pro}\backslash i^{-1}(b)$. $\t{\pro}$ is 
called the blow up of $\pro$ at the point $b$.
 Roughly speaking, we delete the point $b$ from $\pro$ and 
 substitute it by a projective space of dimension $n-1$.

\begin{theo}
\label{monodromy}
Suppose that the pencil $\{H_t\}_{t\in \S}$ is in general position with respect to $X$
and let $b\in\S\backslash C$, then
\begin{enumerate}
\item
For every two Lefschetz vanishing cycles $\delta_0$ and $\delta_1$ in 
$X_b$ there exists a closed path $\lambda$ in $\S\backslash C$ with initial 
and end point $b$ and such that
\[
h_{\lambda}(\delta_0)=\pm \delta_1
\]
where $h_\lambda$ is the monodromy along the path $\lambda$;
\item
 If $H_{n-1}(X\backslash H_a)=0$ for some $a\in \S\backslash C$ and $H_{n-1}(X_b)\not =0$ 
 then for every Lefschetz vanishing cycle $\delta$ in $L_b$, 
 the action of the monodromy group on $\delta$ generates $H_{n-1}(L_b)$.
\end{enumerate}
\end{theo}

\pr 
The first statement and its proof can be found in 7.3.5 of \cite{lam}. But we can give a rather short proof for it as follows:  
\\
Let us consider the pencil $\{H_t\}_{t\in G}$ as the  projective line  $G$ 
in $\check{\pro}$. Let $\delta _0$ and $\delta_1$ vanish along the paths $\lambda_0$ and 
$\lambda_1$ which connect $b$ to critical values $c_0$ and $c_1$ in $G$, 
respectively. The subset $Z\subset \check{X}$ consisting of all points $x$
 such that the line through $x$ and $b$ is not in general position with 
 respect to $\check{X}$ is a proper and algebraic subset of $\check{X}$. Since $\check{X}$ is an
 irreducible variety and $c_0,c_1\in \check{X}\backslash Z$, there is a path 
 $w$ in $\check{X}\backslash Z$ from $c_0$ to $c_1$. Denote by $G_s$ the line 
 through $b$ and $w(s)$. After Blow up at the point $b$ and using the 
 Ehresmann's theorem, we conclude that: 
 \\
There is an isotopy
\[
H:[0,1]\times G\rightarrow \cup_s G_s 
\]
such that  
\begin{enumerate}
\item
 $H(0,.)$ is the identity map; 
 \item 
 for all $s\in [0,1]$, $H(s,.)$ is a $C^{\infty}$ isomorphism between $G$ 
 and $G_s$ which sends points of $\check{X}$ to $\check{X}$ 
 \item
 For all $s\in [0,1]$  $H(s,b)=b$ and $H(s,c_1)=w(s)$
 \end{enumerate}
Let $\lambda_s'=H(s,\lambda_0)$. 
In each Lefschetz pencil $\{H_t\}_{t\in G_s}$ the cycle $\delta_0$ in $X_b$ 
vanishes along the path $\lambda_s'$ in $w(s)$, 
therefore $\delta_0$ vanishes along $\lambda_1'$ in $c_1=w(1)$. 
Consider $\lambda_1$ and $\lambda_1'$ as the paths which start from $b$ and end
 in a point $b_1$ near $c_1$ and put $\lambda=\lambda _1'-\lambda_1$. 
 By uniqueness of the Lefschetz vanishing cycle along a fixed path we
 can see that the path $\lambda$ is the desired path. 
\\
Now let us prove the second part. Unfortunately the above argument is true for the fiber $X_b$ and not  $L_b$. Therefore by Theorem ~\ref{generator} we can
only conclude that the action of the monodromy on a  vanishing cycle generates $H_{n-1}(X_b)$. Since $H_{n-1}(X_b)\not =0$, there is no 
homologous to zero  vanishing cycle in $X_b$ . Let us prove that the intersection matrix $[<\delta_i,\delta_j>]_{r\times r}$ of vanishing cycles is connected i.e., for any two vanishing cycles $\delta$ and $\delta'$ there exists a chain $\delta_{i_1},\delta_{i_2},\ldots,\delta_{i_e}$ of vanishing cycles with the following properties: 
\[
\delta=\delta_{i_1}, \delta'=\delta_{i_e}\ \ \ <\delta_{i_k},\delta_{i_{k+1}}>\not=0,\ \ \ k=1,2,\ldots,e-1
\]
If $\delta'$ is not connected to $\delta$ as above then by Picard-Lefschetz formula $\delta'$ has intersection zero with all cycles obtained by the action of the monodromy on $\delta$. But the action of the monodromy on $\delta$ generates $H_{n-1}(X_b)$. $X_b$ is compact and so $\delta'=0$ in $X_b$ which is a contradiction.  Now using Picard-Lefschetz formula in $L_b$ we see that the action of the monodromy on a vanishing cycle generates any other vanishing cycle in $L_b$. By Theorem ~\ref{generator}, vanishing cycles generate $H_{n-1}(L_b)$ and so the second statement is proved.  $\square$ 
\section{Topology of Integrable Foliations}
\label{top}
\ \ \ \ \ 
In this section we will combine the results of the sections ~\ref{ramsec} and
~\ref{picsec} to generalize Theorem ~\ref{generator} and Theorem
~\ref{monodromy} for the foliation $\F(\INT)$. Note that the first Integral of
$\F$ has the critical fibers $\{F=0\}$ and $\{G=0\}$, if $p>1$ and $q>1$
respectively, which don't appear in the Lefschetz foliation. Homologies are
considered in an arbitrary field except in the mentioned cases.   
\subsection{Integrable Foliations and Lefschetz Pencil}
\ \ \ \ \
Let $\F(\INT)$ be an integrable foliation satisfying the generic conditions
 of Proposition ~\ref{gen}. Put 
\[ 
D_1=\{F=0\},\ D_2=\{G=0\}, L_b=(\integ)^{-1}(b)\backslash \R,\ X_b=L_b\cup \R
\]
Consider the reduced normal crossing divisor $D=D_1+D_2$ 
and the positive integer numbers $q,p$ such that 
\[
deg(F)=qd,\ deg(G)=pd,\ g.c.d.(p,q)=1
\]
It is a well-known fact that the fundamental group of the complement of any smooth hypersurface 
$V$ in $\pro$ is isomorphic to $\Z _{deg(V)}$ and therefore
\[
\pi_1(\pro\backslash D_1)_q=\Z _q\  , \ \pi_1(\pro\backslash D_2)_p=\Z _p
\]
By Theorem ~\ref{ramtheo} there exists a degree $pq$ ramification map 
\begin{equation}
\label{rammap}
\tau :\t{\pro}\rightarrow \pro
\end{equation}
with divisor $D$ and ramification
multiplicities
 $q$ and $p$ in $D_1$ and $D_2$, respectively. We can view the polynomials 
 $F$, $G$ and the coordinates $x_i, i=0,\ldots, n-1$, as meromorphic functions 
 with the pole divisor 
 $H_\infty$, the hyperplane at infinity. 
 We have
\[
div(\t{F})=q(\t{D_1} -d.\t{H}_\infty)
\]
 \[
div(\t{G})=p(\t{D_2} -d.\t{H}_\infty)
\]
therefore $\t{F}^{\frac{1}{q}}$ and $\t{G}^{\frac{1}{p}}$ are well-defined 
meromorphic functions on $\t{\pro}$. Define
\[
j:\tpro\backslash \t{H}_\infty\rightarrow \C^2
\]
\[
j(x)=(\t{F}^{\frac{1}{q}}, \t{G}^{\frac{1}{p}})
\]

The following proposition shows that the
 different sheets of $\t{\pro}$ are due to the different values of 
 $\t{F}^{\frac{1}{q}}$ and $\t{G}^{\frac{1}{p}}$.
\begin{cla}
\label{onetoone}
For any $x\in \pro\backslash H_\infty$ the map $j$ 
takes distinct values in $\tau^{-1}(x)$. (If $x\in H_\infty$ choose 
another hyperplane as the hyperplane at infinity). 
\end{cla}

\pr The set 
\[
S=\{x\in \pro \mid \exists a,b\in \t{\pro} \ s.t.
\ \tau (a)=\tau(b)=x,\ a\not =b,\ j(a)=j(b) \}
\]
 is an open closed subset of $\pro$, because the values of 
 $\t{F}^{\frac{1}{q}}$  ( $\t{G}^{\frac{1}{p}}$) in $\tau ^{-1}(x)$ are the same 
 up to multiplication by some $q$-th ($p$-th) root of the unity. 
 Choosing  normalizing coordinates like in Definition ~\ref{ncd} around the  
 points $a\in D_1\cap D_2$ and $\tau^{-1}(a)$, we have 
\begin{equation}
(x_1, x_2)\stackrel{j}{\leftarrow}(x_1,x_2,y)\stackrel{\tau}{\rightarrow} (x_1^q,x_2^p,y) 
\end{equation}
$\tau$ has the degree $pq$ and so $S$ has not any point near $a$, 
therefore $S$ is empty. $\diamond$
\\

 The foliation $\t{\F}=\tau^{*}(\F)$ in $\t{\pro}$ is also integrable and has
  the first integral $\integr$ with divisor
\[
div(\integr)=\t{D_1}-\t{D_2}
\]
For every $\t{b}\in \S$, let 
\[
\t{L}_{\t{b}}=(\integr)^{-1}(\t{b})\backslash
\t{{\cal R}},\  \t{X}_{\t{b}}=\t{L}_{\t{b}}\cup \t{{\cal R}} ,\  \t{{\cal R}}=\t{D}_1\cap\t{D}_2
\]
The following proposition states the relations between the leaves of $\F$ 
and $\t{\F}$. 
\begin{cla}
\label{a}
The following statements are true:
\begin{enumerate}
\item
$\tau$ maps $\t{\R}$ to $\R$ biholomorphically;
\item
$\tau\mid_{\t{L}_0}:\t{L}_0\rightarrow L_0$
($\tau\mid_{\t{L}_\infty}:\t{L}_\infty\rightarrow L_\infty$) 
is a finite covering map of degree $q$ (repectively $p$);
\item
 For any $c\not =0,\infty$, $\tau$ maps $\t{L_c}$ to $L_{c^{pq}}$ biholomorphically.
\end{enumerate}
\end{cla}

\pr The first and second statements are the results obtained in Proposition 
~\ref{subncd2}.
 For the third it is enough to prove that $\tau\mid_{\t{L_c}}$ is one to one.   
\\
If $\tau(x)=\tau(y)$ and $\integr(x)=\integr(y)$ then 
$\frac{\t{F}^{\frac{1}{q}}(x)}{\t{F}^{\frac{1}{q}}(y)}=
\frac{\t{G}^{\frac{1}{p}}(x)}{\t{G}^{\frac{1}{p}}(y)}$
is a constant which is $p$-th and $q$-th root of the unity, but $g.c.d.(p,q)=1$ 
and so  $\t{F}^{\frac{1}{q}}(x)=\t{F}^{\frac{1}{q}}(y)$ and 
$\t{G}^{\frac{1}{p}}(x)=\t{G}^{\frac{1}{p}}(x)$. By Proposition 
~\ref{onetoone} we conclude that $x=y$. $\diamond$ 
\\
Define 
\[
v:\t{\pro}\rightarrow \pn
\]
\[
v(A)=[\ldots; x_0^{i_0}x_1^{i_1}\cdots x_n^{i_n};\ldots;x_n^d;\t{F}^\frac{1}{q}(A);\t{G}^\frac{1}{p}(A)],\ i_0+\cdots +i_n=d 
\]
$N-2$ is the number of monomials of degree $d$ and the 
variables $x_0,x_1,\ldots,x_n$.
\begin{cla}
\label{embedding}
$v$ is an embedding.
\end{cla}

\pr Consider the following commutative diagram:
\begin{equation}
\begin{array}{ccc}
\label{di1}
\pro & \stackrel{v_d}{\rightarrow} & \C P(N-2) \\
\tau\uparrow &  & i\uparrow  \\
\t{\pro} & \stackrel{v}{\rightarrow} & \pn 
\end{array}
\end{equation}
where $v_d$ is the well-known veronese embedding and $i$ is the projection 
on the first $N-1$ coordinates. 
\\
1. $v$ is one to one;
\\
If $a,b\in \t{\pn}$ and $v(a)=v(b)$ then $v_d(\tau(a))=v_d(\tau(b))$ and 
so $\tau(a)=\tau(b)$ and by Proposition ~\ref{onetoone} we conclude 
that $a=b$.
\\ 
2. $v$ is locally embedding;
\\
For any $a\in \t{\pro}$ choose normalizing coordinates around $a$ and $\tau(a)$. 
For example , if $a\in D_1\cap D_2$ the diagram  ~(\ref{di1}) has the form 
\begin{equation}
\begin{array}{ccccc}
\label{di2}
(x_1^q,x_2^p,y) & \stackrel{v_d}{\rightarrow} & v_d(x_1^q,x_2^p,y) \\
\uparrow &  & \uparrow \\
(x_1,x_2,y) & \stackrel{v}{\rightarrow} & (v_d(x_1^q,x_2^p,y),x_1,x_2)
\end{array}
\end{equation}
we have to prove that the bottom map is an embedding at $0$.
\[
Dv(0)=
	 \left [ \begin{array}{ccc}	
	* & * 	& \frac{\partial v_d}{\partial y} \\
		1 & 0 & 0 \\
		0 & 1 & 0 	
       \end{array} 
	\right ].
\]
$v_d$ is the veronese embedding and so $Dv(0)$ has the maximal rank rank $n$. 
For  other points the proof is similar. $\diamond$       

The foliation $\t{\F}$ is obtained by hyperplane sections of the following Lefschetz pencil
\[
 \{H_t\}_{t\in \S},\ H_t=\{[x;x_{N};x_{N+1}]\in\pn\mid x_{N}=tx_{N+1}\}
\]
$\t{D}_1$ and $\t{D}_2$ intersect each other transversally in $\t{\cal R}=
\t{D}_1\cap \t{D}_2$, and $\integr$ has nondegenerate critical points with
distinct images,  therefore $\{H_t\}_{t\in \S}$ is in general 
position with respect to $X=v(\t{\pro})$.

Now consider the following commutative diagram
\begin{equation}
\label{di3}
\begin{array}{ccc}
\t{\pn}& \stackrel{\tau}{\rightarrow} & \pro \\
\integr\downarrow &  & \integ\downarrow \\
\S & \stackrel{i}{\rightarrow} & \S
\end{array}
\end{equation}
where $i(z)=z^{pq}$. Let $\t{C}$ denote the set of critical values of
$\integr$, then by Proposition ~\ref{a}, we conclude that

\begin{coro}
$\integ$ and $\integr$ are fiber bundle maps over 
$\S\backslash (C\cup A)$ and $\S\backslash \t{C}$, respectively.
\end{coro}

\begin{coro}
\label{kufteh}
Let $b\in \S$ be a regular value of $\integ$. Then for every two Lefschetz vanishing cycles $\delta_1$ and $\delta_2$ in 
$X_b$ there is a
  monodromy $h_\lambda$ such that
\[
h_\lambda (\delta_1)=\pm \delta _2
\]
\end{coro}

\pr Fix a point $\t{b}\in i^{-1}(b)$.
 By diagram ~\ref{di3}, we have the following commutative diagram
\begin{equation}
\label{dimo}
\begin{array}{ccccc}
\pi_1(\S\backslash \t{C},\t{b}) & \times & H_*(\t{X}_{\t{b}}) & 
\rightarrow & H_*(\t{X}_{\t{b}}) \\ 
i_*\downarrow & & \tau_*\downarrow & & \tau_*\downarrow \\
\pi_1(\S\backslash (C\cup A),b) & \times & H_*(X_b) & \rightarrow & 
H_*(X_b)
\end{array}
\end{equation}
$\t{\delta}_i=\tau^*(\delta_i),\ i=1,2$ are two Lefschetz vanishing cycles in
$\t{L}_{\t{b}}$. By Theorem ~\ref{monodromy}, there exists a path $\t{\lambda}\in
\pi_1(\S\backslash \t{C},\t{b})$ such that the related monodromy takes 
$\t{\delta}_1$ to $\pm \t{\delta}_2$. We can assume that this path doesn't pass
through $0$ and $\infty$. Now by the above diagram the path
$i(\t{\lambda})$ is the desired path.$\diamond$

\subsection{More About the  Topology of Integrable Foliation}
\def\intag{\frac{\t{F}}{\t{G}^\frac{q}{p}}}
\ \ \ \ \
Here we want to prove a theorem similar to Theorem ~\ref{ggenerator} for the
foliation $\F(\INT)$.
\\
Let $\tau: \t{\pro}\rightarrow \pro$ be a ramification map with simple 
divisor $D=\{G=0\}$ and multiplicity $p$ at $D$. 
 \[
div(\t{G})=p(\t{D} -d.\t{H}_\infty)
\]
therefore $\t{G}^\frac{1}{p}$ is a well-defined meromorphic function on 
$\t{\pro}$. We denote by $\t{C}$ the set of critical values of $\intag$ in
$\S\backslash\{\infty\}$. Also 
\[
\t{{\cal R}}=\{\t{F}=0\}\cap \{\t{G}=0\}
\] 
The foliation $\t{\F}=\tau^*(\F)$ has the first integral 
$\intag$. Note that $0\in\S$ is no more a critical point of $\intag$. 
Consider the following commutative diagram
\[
\begin{array}{ccc}
\label{di3}
\tpro& \stackrel{\tau}{\rightarrow} & \pro \\
\intag\downarrow &  & \integ\downarrow \\
\S & \stackrel{i}{\rightarrow} & \S
\end{array}
\]
where $i(z)=z^{p}$. Like before we have

\begin{cla}
The following statements are true
\begin{enumerate}
\item
$\intag$ is a fiber bundle map over $\S\backslash (\t{C}\cup\{\infty\})$.
\item
$\tau$ maps $\t{\R}$ ($\t{L}_\infty$) to $\R$ (respectively $L_\infty$) 
biholomorphically;
\item
$\tau\mid_ {\t{L}_0}:\t{L}_0\rightarrow L_0$  
is a finite covering map of degree $q$;
\item
 For any $c\not =0,\infty$, $\tau$ maps $\t{L_c}$ to $L_{c^{p}}$ 
 biholomorphically.
\end{enumerate}
\end{cla}

\begin{theo}
\label{16apr}
If $n=2$ then a distinguished set of 
Lefschetz vanishing cycles related to the critical points in the set
$\t{C}$ generates the first homology group of a regular fiber $\tilde{L}_b$ of $\intag$.
\end{theo}

\pr (n=2) By Theorem ~\ref{generator} it is enough to prove that 
$H_{n-1}(\t{\pro}\backslash \t{D})=0$.
According to Proposition ~\ref{cover},
$\tau_*:H_{n-1}(\tpro \backslash \t{D},\Z)\rightarrow H_{n-1}(\pro\backslash
D,\Z)$ 
is one to one, and we also know  that 
$H_{n-1}(\pro\backslash D, \Z)=\Z_{deg(G)}$, which implies that $H_{n-1}(\pro\backslash D)=0$ in an arbitrary field. These facts imply what we want.
$\diamond$  
\begin{coro}
\label{21agu}
Let $b$ be a regular value of $\integ$ and $\Delta$ be a set of distinguished
Lefschetz vanishing cycles related to the critical points in the set $C$. Let
also $h$ be the monodromy around the critical value $0$. Then the set  
\[
\Delta\cup h(\Delta)\cup\cdots\cup h^{p-1}(\Delta)
\]
generates $H_{n-1}(L_b)$.  
\end{coro}
   
\pr Let $\t{\Delta}$ be a distinguished set of Lefschetz vanishing cycles as in 
Theorem  ~\ref{16apr}. We can see easily that 
$\tau(\t{\Delta})=\Delta\cup h(\Delta)\cup\cdots\cup h^{p-1}(\Delta)$.
$\diamond$
\\

The fiber $L_b$ does not contain the points of $\{F=0\}\cap\{G=0\}$, so this corollary partially claims that the cycle around a point of $\{F=0\}\cap\{G=0\}$ is a rational sum of vanishing cycles. 
In the initial steps of this article my objective was to prove the following 
corollary.
\begin{coro} 
Suppose that $n=2$ and the generic fiber of $\integ$ has genus greater than zero.  Then the action of the monodromy group on a vanishing cycle generates $H_{n-1}(L_b)$.
Let $\omega_1$ be a meromorphic 1-form in the projective space of dimension two whose pole divisor is a union of some fibers of $\integ$. If
\[
  \int _{\delta_t} \omega _1 =0
\]
for a continuous family $\delta_t$ of vanishing cycles, then $\omega_1$ restricted to the closure of each fiber of $\integ$ is exact.
\end{coro}
 We recall that in the above corollary we have assumed  the generic conditions of Proposition ~\ref{gen}.
\\
\pr 
 The first part is a direct consequence of Theorem ~\ref{monodromy} and Proposition ~\ref{a}. For the second part it is enough to prove that 
\[
  \int _{\delta_t} \omega _1 =0
\]
For all 1-cycles in the fibers of $\integ$.$\square$  
\\
Using the ramification map $\tau$, the reader can verify that:
\begin{cla}
\label{khaste}
Let $D_0$ be a small disk around $0$ and $l$ be the straight line which
connects 0 to $b_0$, a point in $\partial D_0$, then
\begin{enumerate}
\item
 $(L_l,L_{b_0})$ is a strong deformation retract  of $(L_{D_0},L_{b_0})$;
\item
There is a $C^\infty$ function $\phi :l\times L_{b}\rightarrow L_{l}$ such that
 $\phi$ is a fiber bundle trivialization on $l\backslash \{0\}$ and the
  restriction of $\phi$ to $\{0\}\times L_{b_0}$, namely $g$, is a finite 
  covering map  of degree $p$ from $L_{b_0}$ to $L_0$ ;
\item
There is a monodromy $h:L_{b_0}\rightarrow L_{b_0}$ around $0$ 
such that for every 
$x\in L_{b_0}$ we have 
\[
g^{-1}(g(x))=\{x,h(x),\cdots ,h^{p-1}(x)\}
\]
in particular $h^p=I$ and $g\circ h=g$.   
\end{enumerate}
\end{cla}
  
 \appendix
\section{Generic Properties}
\label{secgeneric}
\ \ \ \ \
Here we will prove  Proposition ~\ref{gen}. The main tool is the 
transversality theorem which appears both in Algebraic Geometry and Differential 
Topology. We will work in the category of algebraic varieties but
the whole of this
 discussion can be done in the $C^{\infty}$ category of  manifolds.
 \\
  In the sequel by $TX$ we denote the tangent bundle of the variety $X$ 
  and by $(TX)_0$ we denote the image of the zero section of the vector bundle 
  $TX$. For any $x\in X$ we have 
\[
T_{0_x}(TX)=T_{0_x}(TX)_0\oplus T_{0_x}(T_xX)
\]
and so we can define
\[
 d:T_{0_x}(TX)\rightarrow  T_{0_x}(T_xX)
\]
$d$ is the projection on the second coordinate.
 We will essentially use the following transversality theorem in algebraic
  geometry:
\begin{Theo}
\label{transversality}
Let $f:X\rightarrow Z$ and $\pi : X\rightarrow A$ be morphisms 
($C^\infty$ functions) between smooth varieties (resp. $C^\infty$ manifolds) 
and $W$ be a smooth subvariety (resp. submanifold) of $Z$. 
Also assume that $\pi$ is surjective and $f$ is transverse to $W$, then there 
exists an open dense subset $U$ of $A$ such that 
$f\mid _{\pi ^{-1}(\alpha )}$ is transverse to $W$ for 
every $\alpha \in U$.  
\end{Theo}

\pr This theorem is a consequence of Bertini's theorems ( see \cite{sha}). For
more information about the transversality theorem the
reader is referred to \cite{spe} and \cite{abr}.  
\\
Recall that $f:X\rightarrow Z$ is transverse to $W$ if for every $x\in X$ 
with $y=f(x)\in W$, we have $T_{y}W +(T_xf)(T_x(X))=T_y(Z)$. 
This is equivalent to this fact that $f^{-1}(W)$ is empty or is a smooth
 subvariety of $X$ of dimension $dim(X)-dim(Z)+dim(W)$. The following well-known proposition will be used.
\begin{cla}
\label{nondegenerate}
Let $f: X\rightarrow Z $ be a morphism between two smooth varieties and $dim(Z)=1$. Then the critical points of $f$ are nondegenerate, if and only if, $Tf:TX\backslash{(TX)}_0 \rightarrow TZ$ is transverse to ${(TZ)}_0$.
\end{cla}

 Let 
\[
X=\{(F,G,x)\in \PP \times \pro \mid F(x)\not =0 ,\ G(x)\not =0\}
\]
\[
g:X \rightarrow \S,\ g(F,G,x)=\integ (x)=f(x)
\]
and $\tilde{T} X$ be the subvector bundle of $TX$ whose fiber 
$\tilde{T}_{(F,G,x)}X$ is the tangent space of $\{(F,G)\}\times \pro$.
 Let also $\tilde{T}g$ be the restriction of $Tg$ to $\tilde{T} X$ and 
 $\pi :\tilde{T}X\rightarrow \PP$ be the projection on the parameter $(F,G)$.
\begin{cla}
\label{22agu}
For a generic pair $(F,G)$, the critical points of $\integ$ in
$\pro\backslash(\{F=0\}\cup\{G=0\})$ are nondegenerate.
\end{cla}

\pr
According to the transversality theorem and Proposition 
~\ref{nondegenerate} it is enough
 to prove that 
 \[
 \tilde{T}g: \tilde{T}X \rightarrow T\S
 \]
is transverse to $(T\S )_0$.  In a local coordinate around 
$(F,G,x',v)\in \tilde{T}X$ we have
\[
\tilde{T}g(F,G,x',v)=(\integ (x'),D(\integ) (x')(v))
\]
\[
B=d\circ T_{(F,G,x',v)}(\tilde{T}g)(\bar{F},\bar{G},u,w)=
\]
\[ D^2f(x')(v)(u)+Df(x')(w)+pD(f\frac{\bar{F}}{F})(x')(v)-qD(f\frac{\bar{G}}{G})(x')(v)
\]
If $\tilde{T}g$ is not transverse to $(T\S )_0$ at $(x' ,v)$ with 
$v\not = 0$,  then $B=0$ for all $u,w,\bar{F},\bar{G}$. 
Putting $\bar{F}=\bar{G}=0$ we get
\[
Df(x')=0,\ D^2f(x')(v)=0
\]
Let $x'=(x_1',x_2',\ldots x_n')$ and $v'=(v_1,v_2,\ldots v_n)$  
then for all $i=1,2,\ldots ,n$ putting $\bar{F}=x-x'_i,\ \bar{G}=0$,  
we obtain $v_i=0$. This implies that $v=0$ which is a contradiction.$\diamond$
\\

The next step is to prove that generically the images of the critical points 
of $\integ$ are distinct in $\S$. I did not succeed 
to get this generic property by using the transversality theorem, 
therefore I will prove it in the projective space of dimension two, by an elementary arguments in algebraic geometry. The proof in higher dimensions is the same. 
The following lemmas will be used:
\begin{lem}
  \label{linear}
Let $\phi :\C ^n\rightarrow \C ^m $ be a linear map and $A$ be a subvariety of
 $\C ^m$. 
 Then $A\cap Im(\phi)=A_1\cup A_2\cup\cdots$ is the decomposition of 
 $A\cap Im(\phi)$ into irreducible components, if and only if, 
 $\phi ^{-1}(A)=\phi ^{-1}(A_1)\cup \phi ^{-1}(A_2)\cup \cdots$ is the 
 decomposition of  $\phi ^{-1}(A)$ into irreducible components.
\end{lem}

\pr This is due to the fact that we can choose a basis for the vector space 
$\C ^n$ such that
 $\phi :\C ^n =\C ^{n-m'}\times \C ^{m'}\rightarrow Im(\phi)=\C ^{m'}$ 
 is the projection on the second coordinate.$\diamond$
\\

Let 
\[
A'=\{x\in \C ^6 \mid px_4x_2-qx_1x_5=px_4x_3-qx_1x_6=0\}
\]
\[
A'_r=\{x\in \C ^6\mid x_4=x_1=0\}
\]
\[
A'_c=\{x\in \C ^6 \mid px_4x_2-qx_1x_5=px_4x_3-qx_1x_6=x_2 x_6-x_3x_5=0\}
\]
\begin{lem}
\label{ire}
The following statements are true:
\begin{enumerate}
\item
$A'$ has two irreducible components $A'_r$ and $A'_c$;
\item
 $A'\times A'$ has four irreducible components 
 $A'_{ij}=A'_i\times A'_j,\ i,j=r,c$; 
\item 
 For any linear subspace $V$ of $\C^{12}$ of dimension greater than 8,
  $A'\times A'\cap V$ has also four irreducible components 
  $A'_{ij}\cap V,\ i,j=r,c$. 
\end{enumerate}
\end{lem}
In fact, from this lemma we only need to the fact that, for any linear 
subspace $V$ of $\C^{12}$ of dimension greater than 8, 
$A'_{cc}\cap V$ is irreducible. 
\\
\def\pro2{{\Bbb C}P(2)}
Consider an affine open set $\C ^2\subset \pro2$ and let 
\[
0=(0,0),\ 1=(0,1)
\]
Define
\[
A=\{\omega=\INT\mid (F,G)\in \PP\ \&\ \omega   
\hbox{ has singularity at 0 and 1}\}
\]

\begin{lem}
The variety $A$ has exactly four irreducible components 
$A_{rr},\ A_{rc},A_{cr},A_{cc}$. The component $A_{rc}$ contains 
all 1-forms in $A$ which have a radial singularity at $0$ and a center 
singularity at $1$. In the same way other components are defined. 
\end{lem}
\pr
For any $p\in \C ^2$ define the linear map
\[
	\phi _{p}:\PP\rightarrow \C ^6
	,\ \phi _p(F,G)=(F(p),F_x(p),F_y(p),G(p),G_x(p),G_y(p))
\]
where the partial derivatives are considered in the fixed affine coordinate.
Also we define 
\[	
	\phi :\PP\rightarrow \C ^{12},\ \phi =(\phi _{(0,0)},\phi _{(0,1)} )
\]

We can assume that $deg(F)\geq 2$ and $deg(G)\geq 1$. With this hypotheses 
the reader can check that
$dim(Im(\phi))\geq 8$. Now our assertion is the direct consequence of  Lemmas 
~\ref{linear}, ~\ref{ire}.
\\

{\bf Proof of Proposition ~\ref{gen}:} According to Proposition ~\ref{22agu}, it is enough to prove that generically the image of nondegenerate critical points are distinct. Let
\[
S=\{(F,G)\in A_{cc} \mid \integ (0)=\integ (1)\} 
\]
Let $(F,G)\in\PP$ and $\integ$ have $r$ nondegenerate critical points 
$p_1,\cdots ,p_r$. There is an small perturbation $(\bar{F},\bar{G})$ of
 $(F,G)$ such that $\integb $ has r distinct critical values.
  Suppose that this is not true, then we can assume that $\integ$ has maximal 
  number $r'$ of critical values in some neighborhood of $(F,G)$ and $r'<r$. 
  There exist two critical points $p_1,p_2$ of $\integ$ such that 
  $\integ (p_1)=\integ (p_2)$ and for any $(\bar{F},\bar{G})$ near $(F,G)$ 
  with corresponding critical point $\bar{p_1},\bar{p_2}$ near $p_1$ and $p_2$,
   respectively, we have
\[
\integb (\bar{p_1})=\integb(\bar{p_2})
\]
Let $L$ be the linear automorphism of $\pro2$ which sends $0$ and $1$ to
 $\bar{p_1}$ and $\bar{p_2}$, respectively. 
 In some neighborhood $U$ of $(F\circ L,G\circ L )$ in $\PP$ we have 
 $A_{cc}\cap U\subset S\cap U$. Since $A_{cc}$ is an irreducible variety we
 conclude that 
 $A_{cc}\subset S$ which is contradiction because 
\[
(xy^{a-1}+y^a, \frac{q}{p}-a +x+ay)\in A_{cc}\backslash S 
\]
$\diamond$


\smallskip
\rightline{Hossein Movasati}
\rightline{Instituto de Matem\'atica Pura e Aplicada, IMPA}
\rightline{Estrada Dona Castorina, 110, 22460-320}
\rightline{Rio de Janeiro, RJ, Brazil}
\rightline{ E-mail: hossein@impa.br}
\end{document}